\documentclass[11pt]{article}

\usepackage{latexsym}
\usepackage{graphicx}
\usepackage{epsfig}
\usepackage[cp1250]{inputenc}

\newpage

\title{NUMERICAL MODELS FOR THE SIMULATION OF \\
       THE FRACTIONAL-ORDER CONTROL SYSTEMS}
\author{ \v{L}ubom\'{\i}r Dor\v{c}\'ak }
\date{ }

\begin{document}

\newpage\thispagestyle{empty}
\begin{centering}
\includegraphics[scale=0.6]{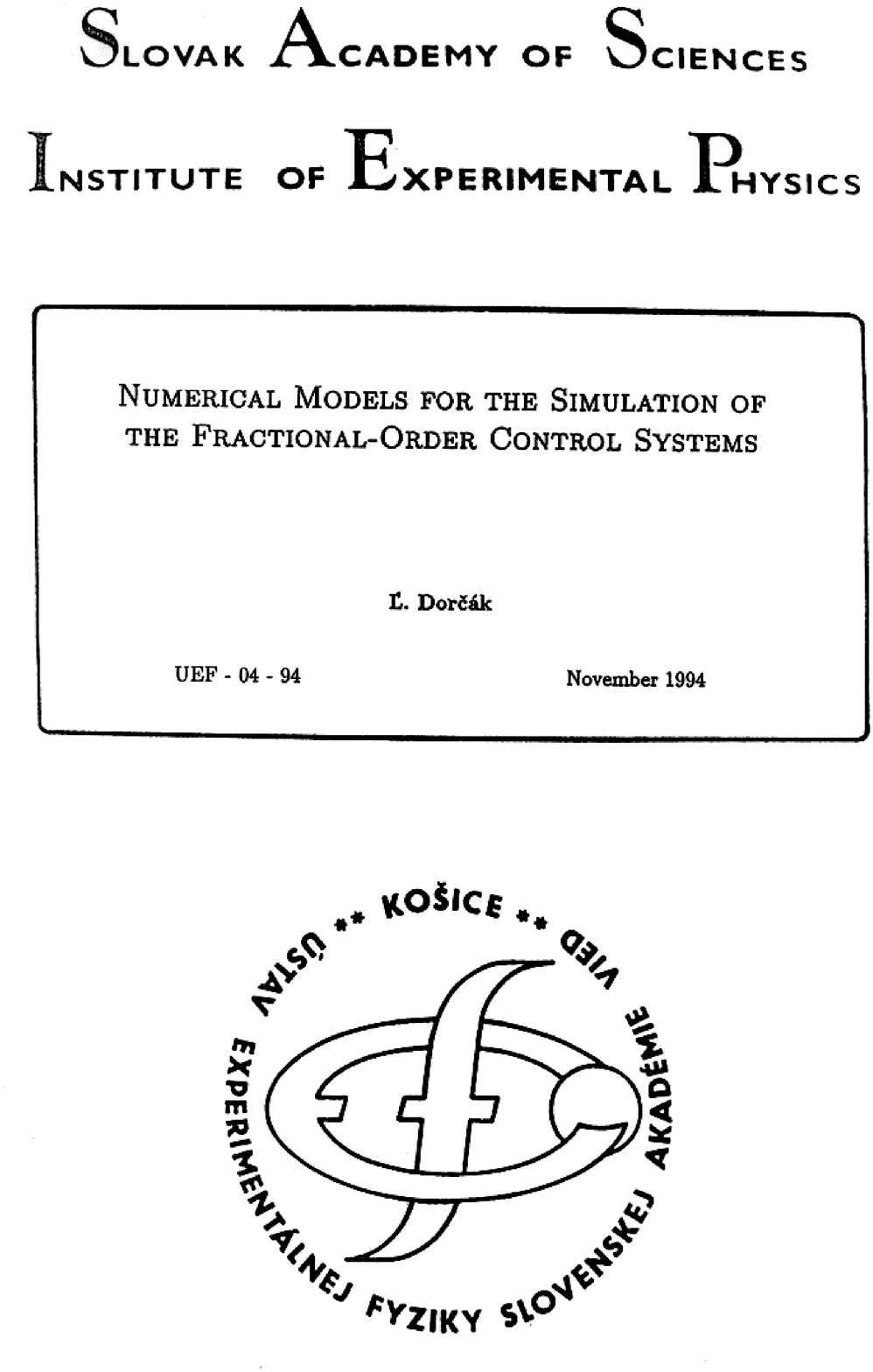}
\end{centering}

\newpage\thispagestyle{empty}
\vspace*{1.5cm}
 \copyright \ 1994, Ing. \v{L}ubom\'{\i}r Dor\v{c}\'ak,~CSc.

\vspace*{20cm}
This publication was typeset by \LaTeX.

\newpage
\setcounter{page}{1}

\vspace*{5cm}
\tableofcontents
\bigskip{\bf References}

\maketitle
\thispagestyle{plain}
\vspace*{0.5cm}
\hspace*{1.7cm} Department of Management and Control Engineering,\\
\hspace*{3.0cm}     BERG Faculty, Technical University of Kosice\\
\hspace*{3.2cm}      Bozeny Nemcovej 3, 042 00 Kosice, Slovakia\\
\hspace*{2cm}  e-mail: Lubomir.Dorcak@tuke.sk,
               phone: (+42155) 6025172 \\

\section{Abstract}

      This contribution deals with the creation of numerical models
for the simulation of the dynamic characteristics of fractional-order
control systems and their comparison with analytical models.
We give the results of the comparison of dynamic properties in
fractional- and integer-order systems with a controller, designed for an
integer-order system as the best approximation to given fractional-order
system.
Other open questions are pointed out, which should be answered in this
area of research.

\section{Introduction}

   The standard control systems used so far were all
considered as integer-order systems, regardless of the reality. In their
analysis and design, the Laplace transform was used heavily for simplicity.
Because of the higher complexity and the absence of adequate mathematical
tools, fractional-order dynamical systems were only treated
marginally in the theory and practice of control systems, e.g.
\cite{L1,L2,L3}.
Their analysis requires familiarity of work with fractional-order
derivatives and integrals
\cite{L4,L5,L6,L9,L10,L11,L12}.

   By removing the restrictions to integer-order systems it is possible
to obtain systems whose properties are a combination of systems of the
closest integer-order, but also intermediate types of systems, which
broadens the class of the systems considerably
\cite{L1}.
  With different fractional-order systems
the notions arise such as weak or strong integrator or differentiator, weak
or strong fractional-type pole, or zero, with interesting contribution to
the dynamics of the system (stability, phase shift etc.), as some properties
are retained, others are eliminated. A fractional-order system combines
some characteristics of systems of the order
 $N$
 and
$ (N+1)$.
  Hence by changing the
order as a real and not only integer value we have more possibilities for an
adjustment of the roots of the characteristic equation according to special
requirements.

      In this contribution we will analyze dynamic properties of systems
in the time domain with an emphasis on the numerical  methods of
simulation of fractional-order systems. We will point out problems of
inadequate approximation of fractional-order systems with
integer-order systems and the differences in dynamic properties of such
systems in closed control systems with controller.

\section{Definition of the fractional-order control system}

   For the definition of  the control system we consider a simple unity
feedback control system shown in Fig.1.
  $G_s(s)$
  denotes the transfer function of
the control system which is either integer-type
 ($G^i_s(s)$)
  or more generally fractional-type
 ($G^f_s(s)$)
 and
 $G_r(s)$
  is the transfer function of the controller, also
either integer-type
 ($G^i_r(s)$)
 or fractional-type
 ($G^f_r(s)$).

\setlength{\unitlength}{1mm}
\begin{picture}(110,42)

\put(25, 2){Figure 1: Simple unity feedback control system}

\put( 6,27){$W(s) \ \ +$}
\put(19,19){$ -$}
\put(29.5,27){$E(s)$}
\put(61.0,27){$U(s)$}
\put(97.0,27){$Y(s)$}
\put( 6,25){\vector(1,0){15.9}}
\put(25,10){\vector(0,1){11.7}}
\put(25,25){\circle{6}}
\put(22.9,27.0){\line(1,-1){4.0}}
\put(22.9,23.0){\line(1, 1){4.0}}
\put(28,25){\vector(1,0){12.2}}
\put(25,10){\line(1,0){72.1}}

\put(40,20){\framebox(20,10)
  [cc]{$G_r(s)$}            }
\put(60,25){\vector(1,0){10.0}}

\put(70,20){\framebox(20,10)
  [cc]{$G_s(s)$}            }
\put(90,25){\vector(1,0){15.9}}
\put(97,10){\line(0,1){15.1}}

\end{picture}

  For later purposes consider a fractional-order controlled system,
which represents our real system, with the transfer function
\begin{equation} \label{r1}
      G^f_s(s) = \frac{1}{a_2s^{\alpha} + a_1s^{\beta} + a_0}
\end{equation}
where
 $\alpha$ and $\beta$
 are in general real ($\alpha > \beta$).
 In the simulation, the coefficient values
 $a_2=0.8,\ a_1=0.5,\ a_0=1,\ \alpha=2.2,\ \beta=0.9$
  were chosen
\cite{L10, L12}.
   To the fractional-type transfer function
 (\ref{r1})
     there
corresponds, in time domain, the fractional-order differential equation
\begin{equation} \label{r2}
      {a_2y^{(\alpha)}(t) + a_1y^{(\beta)} (t) + a_0 y(t)} = u(t)
\end{equation}
with initial conditions
 $y^{(\beta)}(0)=0$ and $y(0)=0$.

\section[Numerical and analytical computation of the \\
         unit step response of a fractional-order system
        ]
        {\protect\parbox[t]{11.5cm}
           {Numerical and analytical computation of the unit
            step response of a fractional-order system
           }
        }

   For the numerical calculation of the
unit step response of the fractional-order system
 (\ref{r2})
we employ, for the approximation of the fractional derivatives in equation
 (\ref{r2}),
the relation (\ref{rr3}) with "short memory" principle, formulated
in  \cite{L11}
\begin{equation} \label{rr3}
   y^{(\alpha)}(t) \approx
   \; _{(t-L)}D^{\alpha}_{t}y(t)\ = \;
   h^{-\alpha} \sum_{j=0}^{N(t)}b_{j}y(t-j h ) \; ,
\end{equation}
where $L$ is "memory length", $h$ is time step,

\begin{displaymath}
N(t) = \min \left\{
           \left[
              \frac{\mbox{$t$}}{\mbox{$h$}}
           \right] , \;
           \left[
              \frac{\mbox{$L$}}{\mbox{$h$}}
           \right]
         \right\} \; ,
\hspace*{1em}
\end{displaymath}
$[z]$ is the integer part $z$,

\begin{equation} \label{r5}
       b_j = (-1)^j {{\alpha} \choose j}
\end{equation}
where
$
        {{\alpha} \choose j}
$
 is binomial coefficient. To calculate
 $b_j$
 it is convenient to use the following recurrent relation
\begin{equation} \label{r6}
       b_0 = 1
         \ , \ \
       b_j = (1- \frac {1+\alpha} {j}) \ b_{j-1}
\end{equation}

It follows from the estimates derived in
 \cite{L11}
that in our case the normed error of such approximation is
\begin{equation}
\delta_{0}= \frac{\left| y^{(\alpha)}(t) - \;
                        _{(t-L)}D^{\alpha}_{t}y(t)\right|}{M} =
      \frac{\mbox{$1$}}
     {\sqrt{\mbox{$L$}} \; \Gamma (\alpha)}\;,
\hspace*{1em}
M=\max_{[0, \; \infty]} |y(t)|
\end{equation}
whence we have the following constraint for the choice of "memory length"
 $L$:
\begin{equation}
     L \geq  \frac{\mbox{$1$}}{\delta_{0}^2 \; \Gamma^2(\alpha)}
\end{equation}
where $\delta _{0}$ is the maximum admissible normalized error
and $\Gamma(\alpha)$ is the {\em Gamma function}.

    By using the relation (\ref{rr3})
    we can approximate the differential equation (\ref{r2})
a different mode.
    Our approximation \cite{L10, L12} of equation
 (\ref{r2})  in discrete time steps
 $t_m\  (m=2,3,...)$   has the following form
\begin{equation} \label{r3}
       a_2h^{-\alpha}\sum_{j=0}^{m} {b_j y_{m-j}} +
       a_1h^{-\beta} \sum_{j=0}^{m} {c_j y_{m-j}} +
       a_0y_{m} = u_m
\end{equation}
or                                                              \\
\begin{equation} \label{r4}
       a_2h^{-\alpha} (b_0 y_{m}+
       \sum_{j=1}^{m} {b_j y_{m-j}}) +
       a_1h^{-\beta}  (c_0 y_{m}+
       \sum_{j=1}^{m} {c_j y_{m-j}}) +
       a_0y_{m} = u_m
\end{equation}

  From the approximation (\ref{r4}) we can derive \cite{L10, L12},
 the following explicit recurrent
relation for the calculation of the values
 $y_m\ (m=2,3,...)$
\begin{equation} \label{r7}
       y_m =
       \frac{u_m -
       a_2h^{-\alpha} \sum\limits_{j=1}^{m} {b_j y_{m-j}} -
       a_1h^{-\beta}  \sum\limits_{j=1}^{m} {c_j y_{m-j}}
            }
            {
       a_2h^{-\alpha} b_0 +
       a_1h^{-\beta}  c_0
     + a_0
            }
\end{equation}
with
\ $y_0 = 0,\ y_1 = 0, \ u_0 = 0 \ $ and $ \ u_m = 1 \ $ for $ \ m=1,2,... $.

    This algorithm does not require iterational calculations, in contrast to
the procedure given in
\cite{L4}.

   For the analytical calculation of the
unit step response of fractional-order systems
 (\ref{r2})
we apply the analytical form of the impulse response of such system
\cite{L5}.
By integrating the impulse response of such system
we obtain the following analytical
form of the
unit step response of such a fractional-order system
\begin{equation} \label{r8}
      y(t) = \frac{1}{a_2}
             \sum_{m=0}^{\infty}
             \frac{(-1)^m}{m!}
             \left( \frac{a_0}{a_2}  \right)^{m}
             t^{\alpha (m+1) }
             E_{\alpha - \beta, \alpha + \beta m+1}^{(m)}
               (-\frac{a_1}{a_2} t^{\alpha - \beta})
\end{equation}
where $E_{\lambda, \mu}(z)$ is the
Mittag-Leffler function in two parameters,
\begin{equation} \label{r9}
E_{\lambda,\mu}^{(k)}(z) \equiv \frac{d^{k}}{dz^{k}}E_{\lambda ,\mu}(z) =
\sum_{j=0}^{\infty} \frac{(j+k)! \,\, z^{j}}
                         {j! \,\, \Gamma (\lambda j + \lambda k + \mu)},
\hspace{3em}
(k = 0, 1, 2, ...)
\end{equation}

    For the calculation of the
  {\em Gamma function}
 we used the algorithm from
\cite{L7}.
The calculation of the
unit step response directly according to
 (\ref{r8})
 is numerically difficult. For the system
 (\ref{r2}),
   at least 20 terms of the series
were needed along with the extended precision of the real variables.
To employ
 (\ref{r8})
  further it is necessary to consider asymptotic methods of
calculations.

\vskip 2.5 mm
\begin{figure}[ht]
\centering
\includegraphics[scale=0.24]{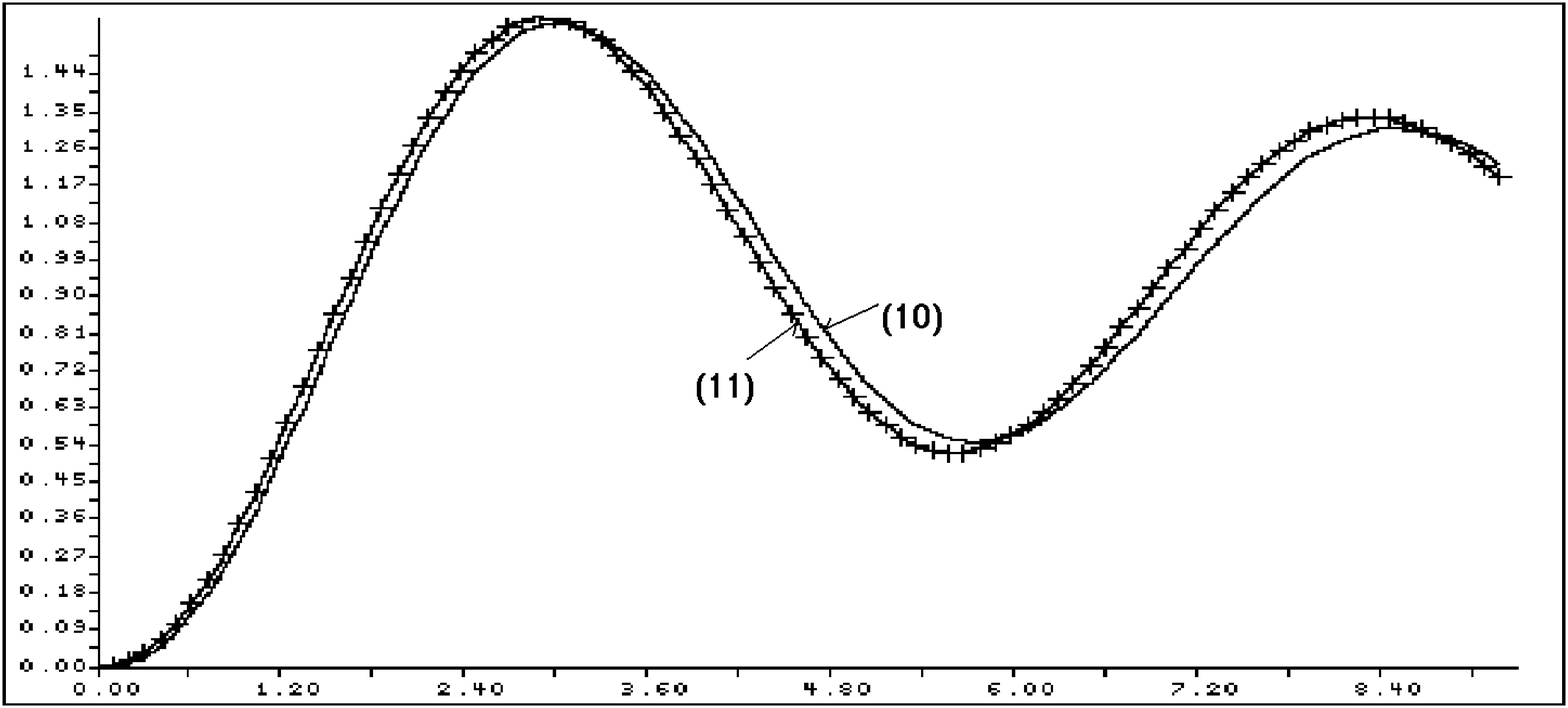}
\vskip 2.5 mm
\centerline{Fig.2 Comparisom of unit sptep responses (10), (11)}
        \label{fig:Fig2}
\end{figure}

\vskip 1.5 mm
   Fig. 2
 shows a graphical comparison of numerically and analytically computed
unit step responses of the system
 (\ref{r2})
 with time step
 $0.1$
 seconds.
The calculations show a good agreement of both methods. The differences
depend on the size of the time step and on the method of approximation
of the value
 $y_1$.

The value $y_1=0$ in  (\ref{r7})  was calculated using the initial
conditions for differential equation  (\ref{r2}) and relation (\ref{rr3})
with $m=1$
\begin{equation} \label{rrr3}
  y_1 = \frac { \ y^{(\beta)}_0 \ -  \ y_0 \ b_1 \ h^{-\beta} }
              { b_0 \ h^{-\beta} }
\end{equation}

Better results (see Fig. 3) we can obtain by approximation the value $y_1$
from initial condition $y_0$ directly with equation
(\ref{r7})
for $m=1$. The advantage of such approximation is more evident for
differential equations
(\ref{r13}), (\ref{r15}) and (\ref{r19}) with derivatives of the function
$u(t)$ (unit step).

\vskip 2.5 mm
\begin{figure}[ht]
\centering
\includegraphics[scale=0.24]{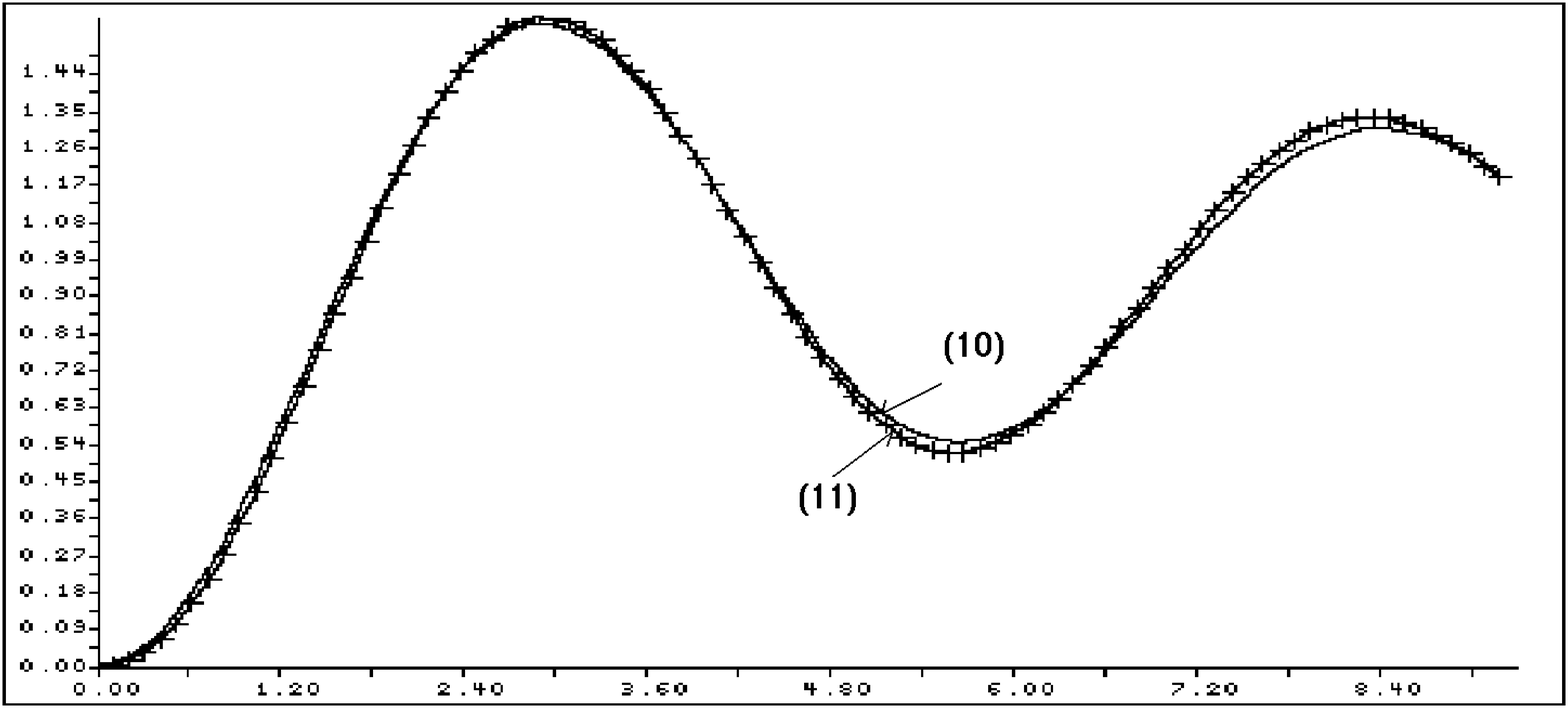}
\vskip 2.5 mm
\centerline{Fig.3 Comparisom of unit sptep responses (10), (11)}
        \label{fig:Fig3}
\end{figure}

\section[Approximation of fractional-order system \\
         with integer-order system
        ]
        {\protect\parbox[t]{10.4cm}
           {Approximation of fractional-order system
            with integer-order system
           }
        }

   In integer-type  linear systems, higher-order systems usually are
approximated - under certain conditions
\cite{L8} -
  with second-order systems for
simplicity. So far, however, fractional-order systems usually have not been
 considered. Regardless of the reality, their approximation is performed
with close integer-type systems, albeit with unfavorable consequences.

   Let us approximate the fractional-type system
(\ref{r2}),
 with the coefficients
as given there, with the integer-type system of second order
\begin{equation} \label{r10}
      {a^i_2 y''(t) + a^i_1 y'(t) + a^i_0 y(t)} = u(t)
\end{equation}

The coefficients
 $a^i_k$
 of the integer-type system take on the values
 $a^i_2=0.7414,\ a^i_1=0.2313,\  a^i_0=1$
under the condition
\begin{equation} \label{r11}
       \sum_{j=0}^{m} {(y^f_j - y^i_j)}^2 \ \  \ \rightarrow \ min
\end{equation}
\vskip 2.5 mm
\begin{figure}[ht]
\centering
\includegraphics[scale=0.24]{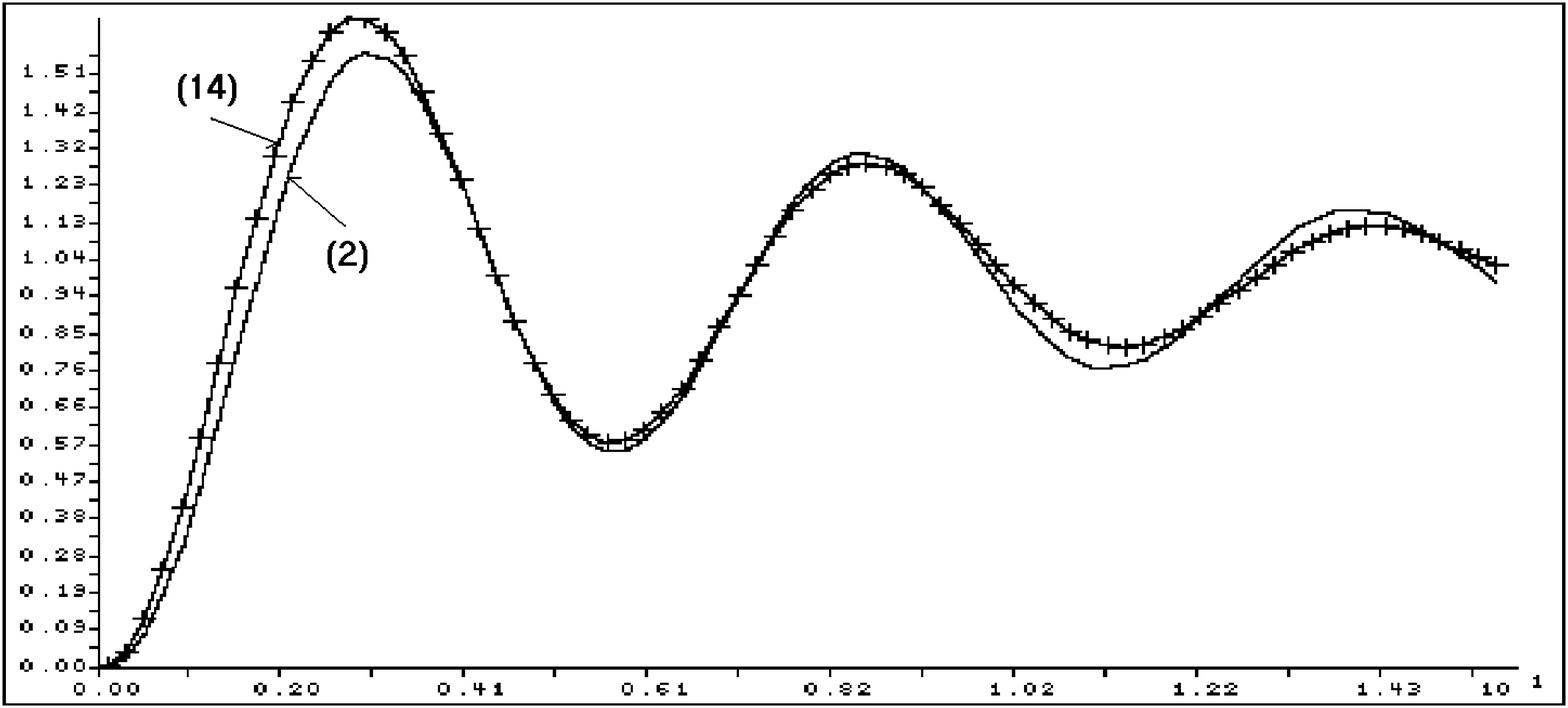}
\vskip 2.5 mm
\centerline{Fig.4 The integer approximation (14) of the system (2)}
        \label{fig:Fig4}
\end{figure}

\vskip 1.5 mm
The comparison of the
unit step responses
of both systems is shown in
Fig.4.
It should be noted that a generally
better approximation was achieved for
$\alpha, \beta \leq 2$.

\section{Design of the regulator}

   For the system
(\ref{r10})
 we choose as an approximation to the system
(\ref{r2})
 the integer-order regulator
\cite{L10, L12}
 with the transfer function
\begin{equation} \label{r12}
       G^i_r(s) = K + T_d \ s
\end{equation}

The regulator will be designed so that a unit step at the input of the
closed regulation system in
 Fig.1
 will induce  at the output an oscillatory
unit step response
with stability measure
 $St=2$
 and damping measure
 $Tl=0.4$.
Then the coefficients for
 (\ref{r12})
 take on the values
 $K=20.5$ and $T_d=2.7343$ .
The differential equation of the integer-order closed regulation system
has the form
\begin{equation} \label{r13}
      {a^i_2 y''(t) + (a^i_1 +T_d) y'(t) + (a^i_0 + K) y(t)} =
      K w(t) + T_d \ w'(t)
\end{equation}

It follows from an analysis of the roots of its characteristic equation that
the requirements for stability measure and damping measure are satisfied.
The permanent regulation deviation is
 $4.6\%$ .
 For the closed regulation system
(\ref{r13})
it is very simple to derive the analytical form of the
unit step response.
In a manner mentioned above it is also possible
to obtain an approximating recurrent relation for the numerical
calculations
\begin{eqnarray} \label{r14}
       y_m =
       \frac{K w_m +
       T_d h^{-1} \sum\limits_{j=0}^{m} {d_j w_{m-j}} +
       a_2 h^{-2} (2y_{m-1} - y_{m-2}) +
       (a_1 + T_d)h^{-1} y_{m-1}
            }
            {
       a_2h^{-2} +
       (a_1 + T_d)h^{-1}
     + (a_0 + K)
            }
\end{eqnarray}
for
 $m=1,2,...$ \ ,
\ \ $y_0 = 0, y_{-1} = 0$,
\ \ $w_0 = 0 \ $ and  $ \ w_m = 1 \ $  for $ \ m=1,2,... $

The
unit step response
is plotted in
 Fig.5.
From the results of numerical simulation it also follows that the regulation
area within
 $10$
seconds was
 $0.71$
 and the permanent regulation deviation is up to
 $5\%.$
\vskip 3.0 mm
\begin{figure}[ht]
\centering
\includegraphics[scale=0.235]{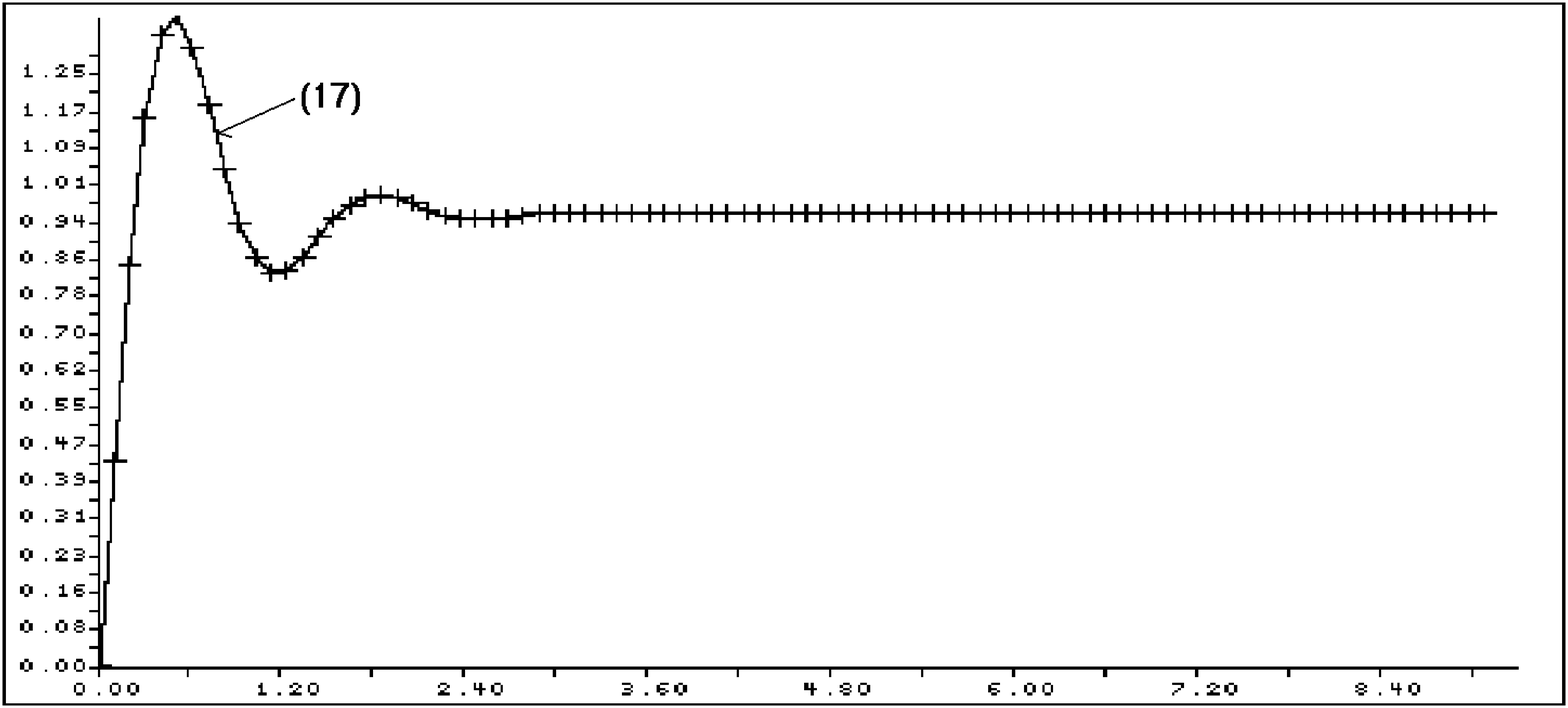}
\vskip 3.0 mm
\centerline{Fig.5 The unit step response (18) of the system (17)}
        \label{fig:Fig5}
\end{figure}

\vskip 3.0 mm
  We now apply the integer-order regulator, designed for the approximated
system (model of the real system), to the original fractional-type
controlled system
(\ref{r2}).
The differential equation of the closed fractional-order regulation system
has the form
\begin{equation} \label{r15}
      {a_2 y^{(\alpha)}(t) + T_d y'(t) + a_1 y^{(\beta)}(t)
      + (a_0 + K) y(t)} = K w(t) + T_d \ w'(t)
\end{equation}

However, the analytical analysis of this equation is not as simple as with
the equation (\ref{r13}).
  We will calculate the unit step response numerically with the aid
of an approximating recurrent relation derived as mentioned above
\begin{eqnarray} \label{r16}
       y_m =
       \frac{K w_m +
       T_d h^{-1}      \sum\limits_{j=0}^{m} {d_j w_{m-j}} -
       a_2 h^{-\alpha} \sum\limits_{j=1}^{m} {b_j y_{m-j}} -
       a_1h^{-\beta}
                      \sum\limits_{j=1}^{m} {c_j y_{m-j}} +
       T_dh^{-1} y_{m-1}
            }
            {
       a_2h^{-\alpha} b_0 +
       a_1h^{-\beta}  c_0 +
       T_dh^{-1}
     + (a_0 + K)
            }
                        \nonumber \\
       \frac{
            }
            {
            }
\end{eqnarray}
for
 $m=1,2,...$ \ ,
\ \ $y_0 = 0$,
\ \ $w_0 = 0 \ $ and  $ \ w_m = 1 \ $  for $ \ m=1,2,... $

To compare the numerical calculations we will also derive an analytical form
of the unit step response of the closed regulating system with a
fractional-order controlled system and an integer-type regulator. We will
obtain the form by integrating the impulse characteristic from
\cite{L5}
 for the
equation of type
 (\ref{r15}),
 whence
\begin{eqnarray} \label{r17}
\hspace{-1em}
      y(t)  =
                 \sum_{m=0}^{\infty}
                 \frac{(-1)^m}{m!}
                 \left( \frac{a_0+K}{a_2}  \right)^{m}
                 \sum_{k=0}^{m} {m \choose k}
                 \left( \frac{a_1}{a_0+K}  \right)^{k}
                        \nonumber \\
\hspace{-1em}
                 ( \frac{ t^{\alpha (m+1) - \beta k} }{b_2}
             E_{\alpha - 1, \alpha + m - \beta k+1}^{(m)}
               (-\frac{T_d}{a_2} t^{\alpha - 1}) +
                         \\
\hspace{-1em}
                 \frac{ t^{\alpha (m+1) - \beta k - 1} }{c_2}
             E_{\alpha - 1, \alpha + m - \beta k}^{(m)}
               (-\frac{T_d}{a_2} t^{\alpha - 1}) ) 
                         \nonumber
\end{eqnarray}
with  \ \ $b_2 = \frac{a_2}{K}$\ \ and \ $c_2 = \frac{a_2}{T_d}$\ .

Fig. 6 shows a graphical comparison of analytically
(\ref{r17})
and numerically
(\ref{r16})
computed unit step responses of the system
(\ref{r15}) with time step 0.1 seconds.
We can see the better results of the numerical calculations by
approximation the value $y_1$ directly with equation
(\ref{r16})
$(a)$ than by approximation with equation
(\ref{rr3}) $(b)$.
\vskip 2.0 mm
\begin{figure}[ht]
\centering
\includegraphics[scale=0.24]{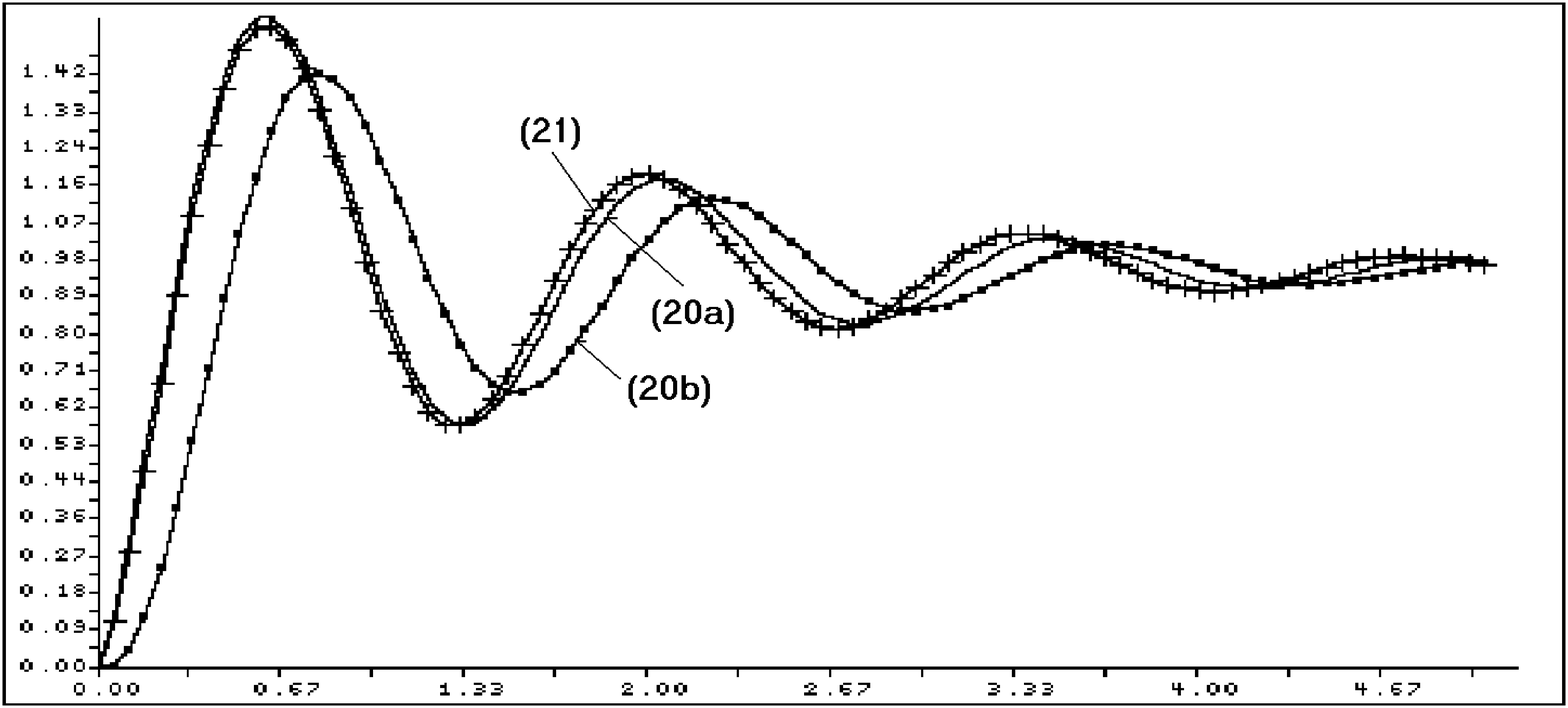}
\vskip 2.5 mm
\centerline{Fig.6 Comparisom of unit sptep responses (20a), (20b) and (21)}
        \label{fig:Fig6}
\end{figure}

\vskip 1.5 mm
A comparison of the unit step responses of the closed regulation system
 with integer-order system
(\ref{r13}),
 with  the fractional-order system
 (\ref{r15})
is shown in
 Fig.7.
\vskip 2.0 mm
\begin{figure}[ht]
\centering
\includegraphics[scale=0.24]{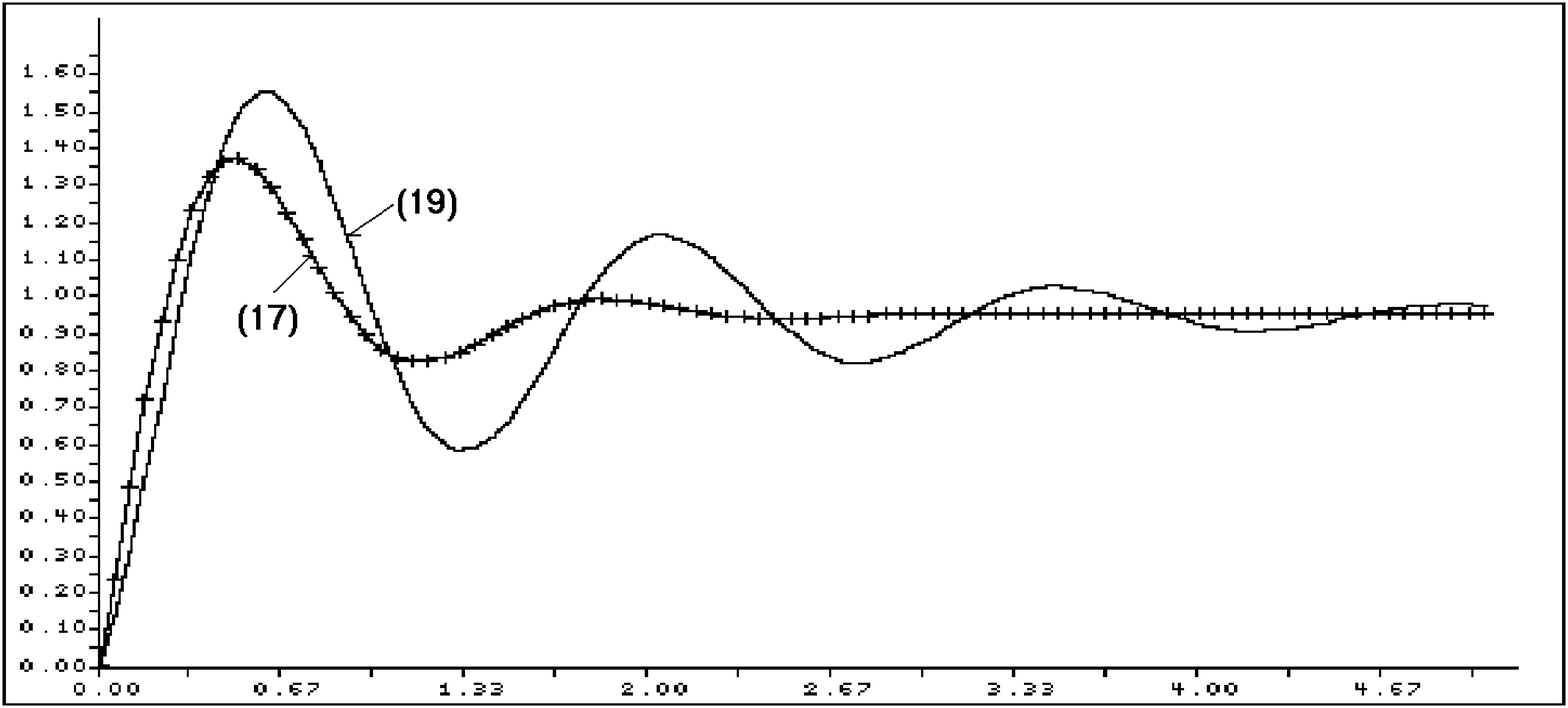}
\vskip 2.5 mm
\centerline{Fig.7 Comparisom of unit sptep responses (17), (19)}
        \label{fig:Fig7}
\end{figure}

\vskip 1.5 mm
It can be seen that the dynamic properties of the closed regulation system
with a fractional-order system and an integer-order regulator designed for
an integer-order system which is an approximation of the fractional-order
system are considerably worse than with the integer-order system and the
integer-order regulator. The regulation area during
 $10$
 seconds is greater by
 $50\%$,
the system stabilizes later and has larger surplus oscillations.
The system
(\ref{r15})
is much more sensitive to changes in parameters.
For example
\cite{L10, L12},
 at the change of
 $T_d$
to value
 $ 1$
the system
(\ref{r15})
is just behind the border of stability
(Fig. 8),
 and with another decrease of
 $ T_d$
it is already unstable, whereas the system
(\ref{r13})
is stable.
\vskip 2.0 mm
\begin{figure}[ht]
\centering
\includegraphics[scale=0.24]{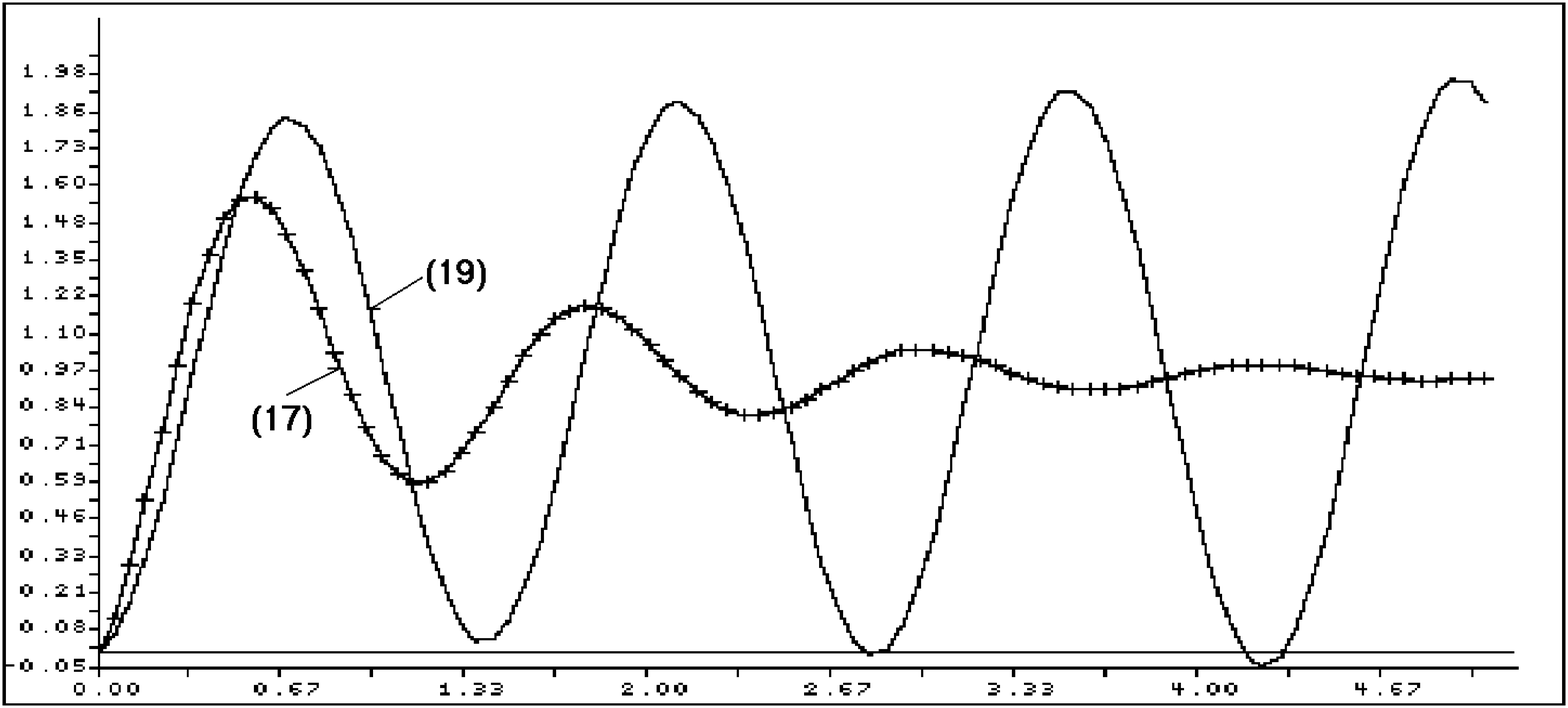}
\vskip 2.5 mm
\centerline{Fig.8 Comparisom of unit sptep responses (17), (19) with
                  changed $T_d$}
        \label{fig:Fig8}
\end{figure}

\vskip 1.5 mm
Hence disregarding the fractional order of the
original system, its approximation by an integer-type system of second order
and an application of a regulator designed for the approximating system to
the original fractional-order system is not adequate for our case
 $ \alpha > 2$.

  For the fractional-order system we now use the fractional-order regulator
\cite{L10, L12}
with the transfer function
\begin{equation} \label{r18}
       G^f_r(s) = K + T_d \ s^{\delta}
\end{equation}

The differential equation of the closed regulation system with a
fractional-order system and a fractional-order regulator has the form
\begin{equation} \label{r19}
  {a_2 y^{(\alpha)}(t) + a_1 y^{(\beta)}(t) + T_d y^{(\delta)}(t)
  + (a_0 + K) y(t)} = K w(t) + T_d \ w^{(\delta)}(t)
\end{equation}

In this case the fractional derivative
 ${\delta}$
of the unit step
 $w^{\delta}(t)$
is no longer equal to zero for $t>0$. The numerical computation of the
unit step response will be done through an approximating recurrent relation
derived as mentioned above
\begin{eqnarray} \label{r20}
       y_m =
       \frac{K w_m +
       T_dh^{-\delta} \sum\limits_{j=0}^{m} {d_j w_{m-j}} -
       a_2h^{-\alpha} \sum\limits_{j=1}^{m} {b_j y_{m-j}} -
            }
            {
       a_2h^{-\alpha} b_0 +
       a_1h^{-\beta}  c_0 +
       T_dh^{-\delta} d_0
     + (a_0+K)
            }
                        \nonumber \\
       \frac{
       a_1h^{-\beta}
                      \sum\limits_{j=1}^{m} {c_j y_{m-j}} -
       T_dh^{-\delta} \sum\limits_{j=1}^{m} {d_j y_{m-j}}
            }
            {
            }
\end{eqnarray}
for
 $m=1,2,...$ \ , \
\ \ $y_0 = 0\ , \ w_0 = 0 $ and $ w_m = 1,\  m = 1, 2, ...$

The analytical form of the unit step response of a closed
regulation system with fractional-order controlled system and
fractional-order regulator can be obtained
\cite{L10, L12}
by dividing equation
(\ref{r19})
into two parts, according to the right-hand side. By taking the derivative
of the order
 ${\delta}$
of the impulse characteristic for the second part,
by subsequent integer-type integration of the impulse characteristic
\cite{L5}
of both parts, and their adding together, we obtain the following resulting
analytical form of the unit step response for equation of the type
 (\ref{r19})
\begin{eqnarray} \label{r21}
\hspace{-1em}
      y(t)  = 
                 \sum_{m=0}^{\infty}
                 \frac{(-1)^m}{m!}
                 \left( \frac{a_0+K}{a_2}  \right)^{m}
                 \sum_{k=0}^{m} {m \choose k}
                 \left( \frac{T_d}{a_0+K}  \right)^{k}
                        \nonumber \\
\hspace{-1em}
                 ( \frac{ t^{\alpha (m+1) - \delta k} }{b_2}
             E_{\alpha - \beta, \alpha + \beta m - \delta k+1}^{(m)}
               (-\frac{a_1}{a_2} t^{\alpha - \beta}) +
                         \\
\hspace{-1em}
                 \frac{ t^{\alpha (m+1) - \delta (k+1)} }{c_2}
             E_{\alpha - \beta, \alpha + \beta m - \delta (k+1)+1}^{(m)}
               (-\frac{a_1}{a_2} t^{\alpha - \beta}) ) 
                         \nonumber
\end{eqnarray}
with  \ \ $b_2 = \frac{a_2}{K}$\ \ and \ $c_2 = \frac{a_2}{T_d}$\ .

   Fig. 9
 shows a graphical comparison of numerically and analytically computed
unit step responses of the system
 (\ref{r19}).

\vskip 3 mm
\begin{figure}[ht]
\centering
\includegraphics[scale=0.24]{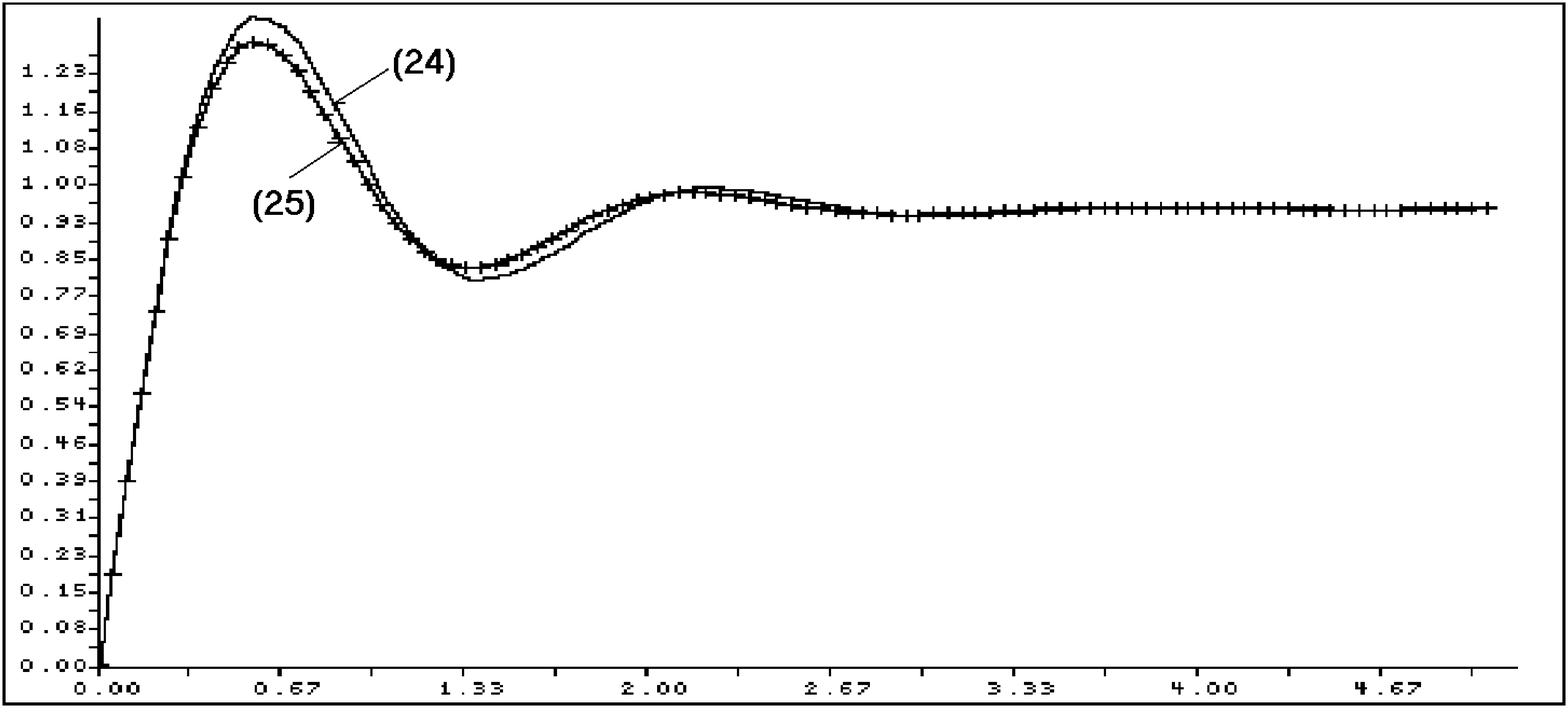}
\vskip 3 mm
\centerline{Fig.9 Comparisom of unit sptep responses (24), (25) of the
                  system (23)}
        \label{fig:Fig9}
\end{figure}

\vskip 2.0 mm
In
 Fig. 10,
a comparison of the unit step responses of the closed regulation
system with a fractional-order regulated system and the integer-order
regulator
(\ref{r15})
and the fractional-order regulator
(\ref{r19}),
 with only two, experimentally found, parameters
$T_d=3.7343$
and
${\delta}=1.15$, is given.

\vskip 3 mm
\begin{figure}[ht]
\centering
\includegraphics[scale=0.24]{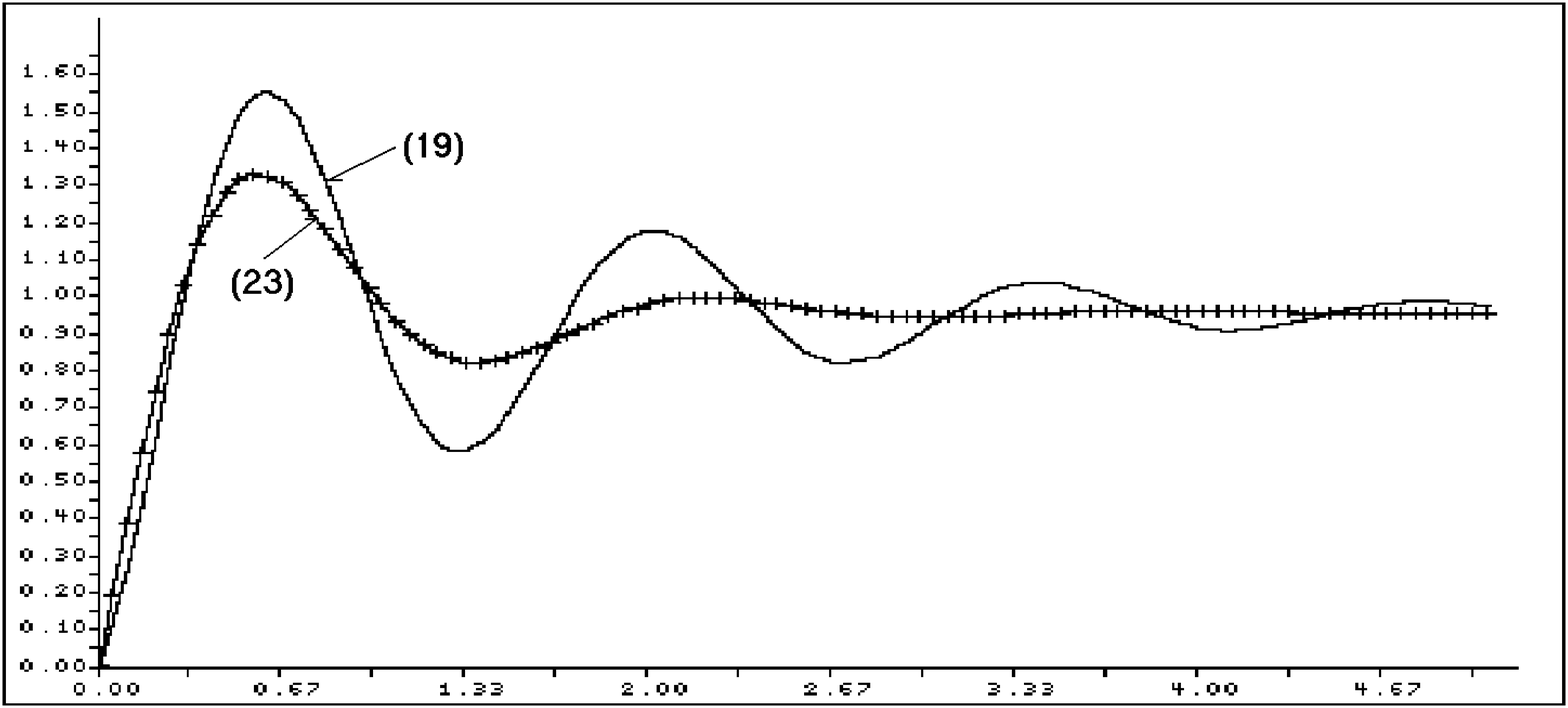}
\vskip 3 mm
\centerline{Fig.10 Comparisom of unit sptep responses of the
                   systems (19) and (23)}
        \label{fig:Fig10}
\end{figure}

\vskip 2.5 mm
It follows from the simulations
(Fig. 10)
that the dynamic properties of the system
(\ref{r15})
with an integer-order regulator (Fig.7)
improved with the use of a fractional-order regulator
(\ref{r19}).
Therefore in the following it is necesary to deal with
the methods of synthesizing the structure and parameters of
fractional-order regulators.

\section{Conclusion}

   In the contribution we gave relations for simulation of dynamic
properties of fractional-order systems, via numerical as well as
analytical methods.
Also presented were the results of some simulations.
The calculations showed a good agreement of the new numerical method and
analytical method of simulation of the fractional-order control systems.

It was shown that the fractional-order system
(\ref{r2})
 can be well approximated by the second-order system
(\ref{r10})
 for
 $ \alpha, \beta \leq 2$.
 The integer-order regulator designed for the integer-order system
(\ref{r10})
can in this case be applied with good results also to the
fractional-order system. A worse approximation is achieved with
 $ \alpha$
 or
 $ \beta > 2$.
An application of the integer-order regulator, designed for the
integer-order system as an approximation to such a fractional-order
system is inadequate and with a change of system or regulator parameters
can lead to system instability.

  In the future it is desirable to deal
with the asymptotic methods
in the analytical method of computation of the dynamic
characteristics of fractional-order systems.
 Further, it is necessary to deal
with the methods of identification of the fractional order and of the
parameters of such systems. The most important task will be the elaboration
of methods of synthesis of the structures and parameters of regulators
for such types of systems.

\newpage\thispagestyle{empty}

\vspace*{17cm}
\noindent
Názov: Numerical \hspace{-0.1mm}Models \hspace{-0.1mm}for
       \hspace{-0.1mm}the \hspace{-0.1mm}Simulation
       \hspace{-0.1mm}of \hspace{-0.1mm}the
       \hspace{-0.1mm}Fractional-Order \hspace{-0.1mm}Control Systems \\
Autor:  Ing. \v{L}ubom\'{\i}r Dor\v{c}\'ak, CSc. \\
Zodp. redaktor: RNDr. P. Samuely, CSc.\\
Vydavate\v{l}: \'Ustav experiment\'alnej fyziky SAV, Ko\v{s}ice \\
Redakcia: \'UEF SAV, Watsonova 3, 04001 Ko\v{s}ice, Slovensk\'a republika \\
Po\v{c}et str\'an: 12 \\
N\'aklad: 120  \\
Rok vydania: 1994 \\
Tlač: OLYMPIA s.r.o.,  M\'anesova 23, 040 01 Ko\v{s}ice

\end{document}